\documentclass[12pt]{article}
\usepackage{amsmath, amsfonts, amsthm, amssymb}
\usepackage{graphicx}
\usepackage{float}
\usepackage{verbatim}% for comment
\usepackage{caption}
\usepackage{subcaption}
\usepackage{color}

\numberwithin{equation}{section}

\newtheorem{thm}{Theorem}[section]
\newtheorem{lma}[thm]{Lemma}

\newtheorem{prop}[thm]{Proposition}

\renewcommand{\ge}{\geqslant}
\renewcommand{\le}{\leqslant}
\renewcommand{\geq}{\geqslant}
\renewcommand{\leq}{\leqslant}
\renewcommand{\H}{\text{H}}

\allowdisplaybreaks

\title{Assouad type dimensions for self-affine sponges with a weak coordinate ordering condition}

\author{Douglas C. Howroyd\\ \\
\emph{School of Mathematics and Statistics,}\\ \emph{ The University of St Andrews,}\\ \emph{St Andrews,} \\ \emph{ KY16 9SS,}\\ \emph{ Scotland}\\ \\
\emph{E-mail contact}:\\  dch8@st-andrews.ac.uk }

\begin{document}

\maketitle

\begin{abstract}
Recently self-affine sponges have been shown to be interesting counter-examples to several previously open problems. One class of recently studied sponges are \emph{Bedford-McMullen sponges with a weak coordinate ordering condition}, that is, sponges with several coordinates having the same contraction ratio. The Assouad type dimensions of such sets cannot be calculated using the same formula as the regular Bedford-McMullen sponges. We calculate the Assouad type dimensions for such sponges and the more general Lalley-Gatzouras sponges with a weak coordinate ordering condition, discussing some of their more subtle details along the way.   \\


\end{abstract}

\section{Introduction} \label{intro}

In the recent paper \cite{fraser-howroyd}, a new phenomenon was noted in the study of self-affine sponges which only occurs in the higher dimensional case (ambient spatial dimensions greater than or equal to 3). The techniques used to calculate the Assouad dimension of Bedford-McMullen sponges in that paper fail for what we will call \emph{sponges with a weak coordinate ordering condition}, where the condition that each coordinate direction is divided into a strictly different number of pieces $n_1 <\ldots < n_d$ is relaxed to $n_1 \le \ldots \le n_d$, which stops us from having a strict ordering of the coordinates. By doing this simple generalisation we will observe a dimension drop which only occurs for the Assouad and lower dimensions. We will also calculate the Assouad and lower dimensions for Lalley-Gatzouras sponges satisfying a weak coordinate ordering condition. This paper will follow the proof in \cite{fraser-howroyd} and, as such, an understanding of that paper is required for the more technical parts of this one.

\subsection{Assouad dimension and lower dimension} \label{dimension}

Many definitions of dimension exist, each with different interesting properties that we would like to study. One such definition is the Assouad dimension which has seen much activity in the past few years, for example \cite{dfsu, fhor, fraser-jordan, mackay}. We will assume throughout that $F \subseteq \mathbb{R}^d$ is a non-empty compact bounded set.  Generally the Assouad dimension of $F$ is defined by 
\begin{multline*}
\dim_{\text{A}} F = \inf \Bigg\{ s \geq 0 \, \,  : \, \exists \text{ constants } C,\rho > 0 \text{ such that, for all } \, \, 0< r< R \leq \rho,\\ \text{ we have  }\sup_{x\in F} N_r (B(x,R)\cap F) \leq C\left(\frac{R}{r}\right)^{s} \Bigg\}\end{multline*}
where $B(x,R)$ means the open ball of radius $R$ and centre $x$ and $N_r(E)$ is the smallest number of open sets in $\mathbb{R}^d$ with diameter less than or equal to $r$ required to cover a bounded set $E$.

When one is considering one of the box dimensions (upper or lower box), it is natural to consider the other at the same time. Similarly the Assouad dimension has a natural dual that we will call the lower dimension $\dim_\text{L} F$, following Bylund and Gudayol \cite{bylund}. Other names do exist for the lower dimension but they all use equivalent definitions. 
\begin{multline*}
\dim_{\text{L}} F = \sup \Bigg\{ s \geq 0 \, \,  : \, \exists \text{ constants }C, \rho > 0 \text{ such that, for all } \, \, 0< r< R \leq \rho,\\ \text{ we have  }\inf_{x\in F} N_r (B(x,R)\cap F) \geq C\left(\frac{R}{r}\right)^{s} \Bigg\}.
\end{multline*}

For more information on some of the properties of these dimensions see \cite{robinson, luk, fraser} but generally the Assouad dimension should be considered as giving us information on the `densest' part of our set whilst the lower dimension does the same for the `thinnest' part, this will become clear in our results. Some of the more common dimensions that exist are the Hausdorff and upper and lower box dimensions, which we denote by $\dim_\H$,  $\overline{\dim}_\text{B}$ and $\underline{\dim}_\text{B}$ respectively, and we refer the reader to \cite[Chapters 2--3]{falconer} for their definitions and basic properties. When the upper and lower box dimension coincide we will simply call the common value the box dimension; this is the case in all sets that are considered in this paper. We will often compare our results to the analogous Hausdorff and box dimension results to highlight interesting properties.  For any compact set $F$, we  generally have 
\[
\dim_\text{L} F  \ \leq \  \dim_\H F\ \leq \  \underline{\dim}_\text{B} F  \ \leq \  \overline{\dim}_\text{B} F \ \leq \  \dim_\text{A} F.
\]

\subsection{Bedford-McMullen sponges with a weak coordinate ordering condition} \label{BMintro}

We will now provide the definition of Bedford-McMullen carpets and sponges, some of the most classical self-affine fractals that have been studied. Much work has been spent understanding the dimension theoretic properties of such sets starting from the original works of Bedford and McMullen \cite{bedford, mcmullen}, followed by \cite{lalley-gatz, baranski, kenyonperes, fengaffine, fraser_box}, to the more recent studies of the Assouad dimension in \cite{mackay, fraser, fraser-howroyd, dfsu, fraser-jordan}. The recent paper \cite{das-simmons} found an interesting application of these results to answer a long existing open problem in dimension theory on the existence of an ergodic invariant measure of full Hausdorff dimension for every expanding repeller, see \cite{kenyonperes} for the full statement. They did this by showing the existence of self-affine sponges where any ergodic invariant measure supported on the sponge has a strictly lower Hausdorff dimension than the sponge itself. Much of the work in this area was aimed at showing that there always existed such a measure, however sponges are simple enough to study yet with complicated enough structure for their result to be possible.

The following notation is a mixture of Olsen's \cite{sponges} work with some additional parts introduced in \cite{fraser-howroyd} modified to solve our problem.  \emph{Bedford--McMullen sponges} are defined as follows. Let $d\in \mathbb{N}$ and, for all $l=1, \ldots, d$, choose $n_l \in \mathbb{N}$ such that $1 < n_1 \le n_2 \le \cdots \le n_d$. When these integers are all equal our sponge is simply a strictly self-similar set and the dimension is equal to the similarity dimension since our sponge satisfies the Open Set Condition, see \cite{olsenassouad}. Our formula will actually work for this case as well, something that the previous formula in \cite{fraser-howroyd} was not able to do. Let $\mathcal{I}_l=\left\{0,\ldots, n_l-1 \right\}$ and $\mathcal{I}=\prod_{l=1}^{d}\mathcal{I}_l$ and consider a fixed digit set $D \subseteq \mathcal{I}$ with at least two elements.  For $\textbf{i} =\left( i_1,\ldots, i_d\right)\in D$ we define the affine contraction $S_\textbf{i}\colon [0,1]^d \rightarrow [0,1]^d$ by
\[
S_{\textbf{i}}\left(x_1, \ldots, x_d \right) = \left( \frac{x_1+i_1}{n_1},\ldots,\frac{x_d+i_d}{n_d} \right) .
\]
Thanks to a well known theorem of Hutchinson $\cite{hutch}$, we know that there exists a unique non-empty compact set $K \subseteq [0,1]^d$ satisfying 
\[
K=\bigcup_{\textbf{i}\in D}S_{\textbf{i}}(K) 
\]
called the attractor of the iterated function system (IFS) $\ \left\{ S_{\textbf{i}} \right\}_{\textbf{i}\in D}$. The attractor $K$ is called a \emph{self-affine sponge} when constructed with these contractions and when $d=2$ we call the resulting set a \emph{carpet}. The $2$ dimensional case was first considered by Bedford and McMullen separately in \cite{bedford, mcmullen} and then Kenyon and Peres \cite{kenyonperes} generalised to the $d>2$ dimensional case. Without loss of generality we shall assume that $K$ does not lie in a hyperplane. In such cases, we simply restrict our attention to the minimal lower hyperplane containing $K$ and consider it as a self-affine sponge in this space.

One advantage of this construction is the ability to link our sponge $K$ and the symbolic space $D^{\mathbb{N}}$, the set of all infinite words over the symbols in $D$ equipped with the product topology generated by the cylinders $[\mathbf{i}_1,\ldots, \mathbf{i}_n]$ corresponding to all finite words over $D$. The function $\tau :D^{\mathbb{N}} \rightarrow [0,1]^d$ defined below is the key to this property:
\[
\tau(\omega)=\bigcap_{n \in \mathbb{N}} S_{\omega\vert n}([0,1]^d)
\]
where $\omega= (\textbf{i}_1, \textbf{i}_2,\ldots)\in D^{\mathbb{N}}$, $\omega \vert n=\left( \textbf{i}_1, \ldots , \textbf{i}_n \right) \in D^n$, $\textbf{i}_j=(i_{j,1},\ldots,i_{j,d})$ for any $j\in \mathbb{N}$, and $S_{\omega\vert n} = S_{\left( \textbf{i}_1,  \ldots ,\textbf{i}_n  \right) } =  S_{ \textbf{i}_1} \circ \cdots  \circ S_{ \textbf{i}_n} $. Thus 
\[
\tau(D^{\mathbb{N}})=K.
\]

The Assouad dimension of Bedford-McMullen sponges was calculated recently in \cite{fraser-howroyd} as long as all of the inequalities in $1<n_1< \cdots < n_d$ are strict. An analogous assumption was assumed in \cite{dfsu} when calculating the Assouad and lower dimensions of Lalley-Gatzouras sponges, a generalisation of the Bedford-McMullen sponges, which indicates that this issue might arise in further generalizations of self-affine sponges. This problem is unique to the Assouad and lower dimensions in the higher dimensional cases ($d>2$) and does not appear in the 2-dimensional case or when considering the Hausdorff or box dimensions of sponges. We call these problematic sponges \emph{Bedford-McMullen sponges with a weak coordinate ordering condition}. Let $1<n_1\le\ldots\le n_d$ with $d>2$ such that there is at least one $1\le a < d$ with $n_a=n_{a+1}$. Using the same method as for the general Bedford-McMullen sponges we obtain an attractor $K$ which we define as a Bedford-McMullen sponge with a weak coordinate ordering condition. The main idea will be to consider the equal coordinates as forming a sort of `strictly self-similar set inside our self-affine sponge' and then employ methods that would generally be used in the self-similar setting for these problematic coordinates whilst continuing the original procedures for the ones with strict inequalities. 

To simplify notation we define $\mathcal{J}_1=\mathcal{I}_1 \times \mathcal{I}_2 \times\cdots \times \mathcal{I}_{a_1}$ such that $1<n_1=\cdots=n_{a_1}<n_{a_1+1}$ then $\mathcal{J}_{2}=\mathcal{I}_{a_1+1} \times \cdots \times \mathcal{I}_{a_1+a_2}$, again with $n_{a_1+1}=\cdots=n_{a_1+a_2}<n_{a_1+a_2+1}$ and we continue by induction on the dimension until we obtain $\mathcal{I}=\prod_{l=1}^{d^*} \mathcal{J}_l$ where $d^*$ is the number of distinct integers from the list $n_1,\ldots, n_d$. So $\mathcal{J}_{l}$ is a collection of $a_l$ dimensions for any $l$ and we define $n_l^*=n_{a_1+\cdots+a_l}$, this is the contraction associated to all the coordinates in the $l^\text{th}$ collection. This groups the `self-similar' parts together and allows us to think of our sponge not as a $d$-dimensional self-affine set without a strict coordinate ordering but a `$d^{*}$-dimensional self-affine set with a strict coordinate ordering' with some `dimensions' actually being several combined coordinates.

A \emph{pre-fractal} of an attractor $F$ associated to an IFS as defined above is the set defined by the application of all possible combinations of functions in our IFS to an initial set a certain number of times (the number of times is called the level of the pre-fractal). As the level tends to infinity, the pre-fractals will converge to our attractor in the Hausdorff metric, for any initial set satisfying some simple conditions, e.g. $\left[0,1\right]^d$. For a more detailed explanation see \cite[page 126]{falconer}. The $n$th pre-fractal of a sponge $K$ is
\[
\bigcup_{\left( \textbf{i}_1,  \ldots ,\textbf{i}_n  \right)\in D^{n} } S_{\left( \textbf{i}_1,  \ldots ,\textbf{i}_n  \right) }([0,1]^d),
\]
which is just a collection of $\lvert D \rvert^n$ rectangles, later we will work with objects that are comparable to the pre-fractals. Below is a figure of the first pre-fractals for a given Bedford-McMullen sponge, we can already see how they will converge to our attractor $K$.    

\begin{figure}[htbp]
\centering
\begin{subfigure}{.5\textwidth}
  \centering
  \includegraphics[width=0.9\linewidth]{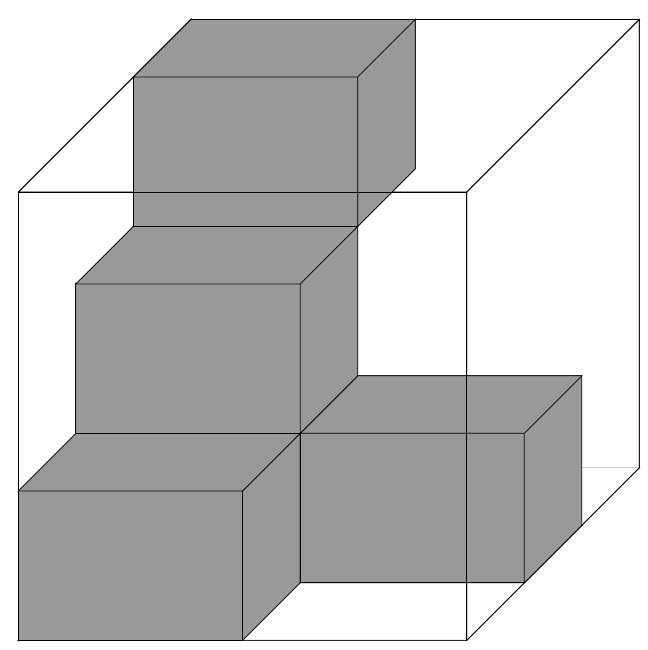}
\end{subfigure}%
\begin{subfigure}{.5\textwidth}
  \centering
  \includegraphics[width=0.9\linewidth]{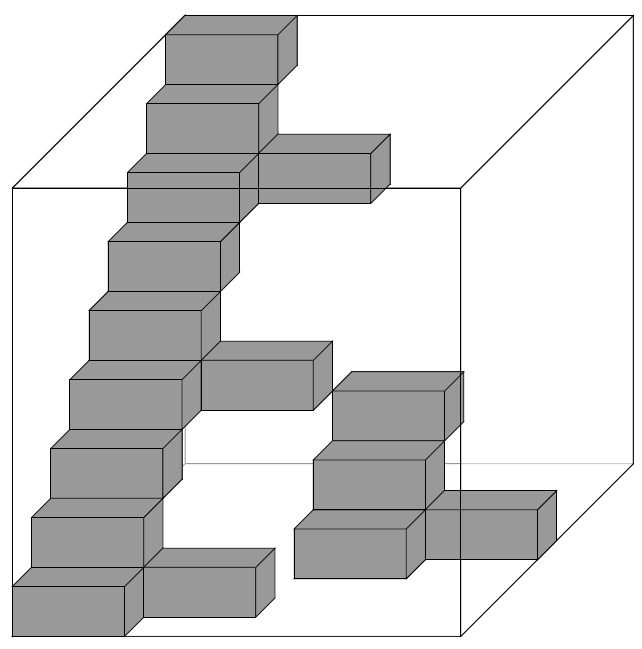}
\end{subfigure}
\caption{The first and second levels in the construction of a specific self-affine sponge satisfying a weak ordering condition in $\mathbb{R}^3$ where $n_1=2$, $n_2=3$, $n_3=3$ and $D = \{(0,0,0), (0,1,1), (0,2,2), (1,0,1) \}$.}
\label{fig:part-affine}
\end{figure}

The technique used to obtain the upper bound for the Assouad dimension will require an understanding of \emph{Bernoulli measures} supported on the set, thankfully iterated function systems provide a nice way of constructing such measures by defining a Borel probability measure on the symbol space and then using the push-forward $\tau$ to acquire a Bernoulli measure on our set. To do this we associate a probability vector $\ \left\{ p_{\textbf{i}} \right\}_{\textbf{i}\in D}$ with $D$ and let $\tilde{\mu}=\prod_{\mathbb{N}}\left( \sum_{\textbf{i}\in D}p_{\textbf{i}}\delta_{\textbf{i}} \right)$  be the natural Borel product probability measure on $D^{\mathbb{N}}$, where $\delta_{\textbf{i}}$ is the Dirac measure on $D$ concentrated at $\textbf{i}$.  Finally, the measure 
\[
\mu(A)=\tilde{\mu}\circ \tau^{-1}(A),
\]
for a Borel set $A\subseteq K$, is our Bernoulli measure supported on $K$. For those not familiar with this technique, one could simply consider the Borel probability measure on the symbol space whilst remembering that cylinders and pre-fractals are closely related. We will later define a  Bernoulli measure which will be used in the proofs to obtain an upper bound on the Assouad dimension. This measure will reflect both the strict inequalities and the equalities of the $n_i$. For now we simply introduce the notation needed for the construction of the measure and the statement of our result.

For any $l=1,2,\ldots,d^*$ we define $\pi_l \, : \, D \rightarrow \prod_{k=1}^l \mathcal{J}_k$ to be the projection onto the first $l$ `clusters' of coordinates so say $\pi_l (i_1,\ldots, i_d)=(i_1,\ldots, i_{a_1+\cdots+a_l})$, let $D_l=\pi_l (D)$ and $N=\#(\pi_1 D)$. Then, for $l=1, \dots, d^*-1$ and $(i_1, \dots , i_{a_1+\cdots+a_l}) \in D_l$ let
\begin{align*}
N(i_1, \ldots, i_{a_1+\cdots+a_l})=& \# \{ \left(i_{a_1+\cdots+a_l+1}, \ldots, i_{a_1+\cdots+a_l+a_{l+1}}  \right) \in \mathcal{J}_{l+1} \\
&\qquad \qquad: (i_1, \ldots, i_{a_1+\cdots+a_l+a_{l+1}}) \in D_{l+1} \}
\end{align*}
 be the number of possible ways to choose the next $a_{l+1}$ digits of $(i_1, \dots , i_{a_1+\cdots+a_l})$. Thus $N(i_1, \ldots, i_l)$ is an integer between 1 and $\left(n_{l}^*\right)^{a_{l+1}}$, inclusive.

\subsection{Lalley-Gatzouras sponges with a weak coordinate ordering condition} \label{LGintro}

One can generalise the construction of Bedford-McMullen sponges to the sets first considered by Lalley and Gatzouras \cite{lalley-gatz}, who calculated the Hausdorff and box dimensions of Lalley-Gatzouras carpets. Mackay \cite{mackay} then calculated the Assouad dimension of these carpets. Following Mackay's notation we let $d\in \mathbb{N}$ and choose a fixed digit set $D\subset \mathcal{I} = \prod_{l=1}^d \left\{ 0,\ldots,n_l-1\right\}= \prod_{l=1}^d \mathcal{I}_l$ where $1<n_1\le n_2 \le \cdots \le n_d$ are integers, note these integers are not always going to be related to our contraction ratios, we simply need them to construct our symbol space. Given an $\textbf{i}=(i_1,\ldots,i_d) \in D$ we define 
\[
D_{i_1,\ldots,i_{l-1}} = \left\{(j_1,\ldots,j_{l-1},j_l) \in \pi_l(D)\colon (i_1,\ldots, i_{l-1})=(j_1,\ldots,j_{l-1}) \right\}
\] 
for any $l=2,\ldots, d$ , where $\pi_l$ is the projection function onto $\prod_{k=1}^l \mathcal{I}_k$. Then for each $\textbf{i}=(i_1,\ldots,i_d) \in D$ we associate contractions $c_{i_1,\ldots,i_l}$ and translations $t_{i_1,\ldots, i_l}$ for all $l=1,\ldots, d$ such that the following conditions hold:
\begin{itemize}
  \item $1>c_{i_1} \ge c_{i_1,i_2} \ge \cdots \ge c_{i_1,\ldots, i_d}>0$,
  \item if $\textbf{i},\textbf{j}\in D $ such that  $(i_{1}\ldots,i_{l}) = (j_1,\ldots j_l)$ for some $l$, then  $c_{i_1,\ldots,i_l}=c_{j_1\ldots,j_l}$,
  \item $\sum_{i\in \pi_1(D)} c_i \le 1$,
  \item for any $l=2,\ldots,d$ and $(i_1,\ldots,i_{l-1}) \in \pi_{l-1}(D)$, $\sum_{(j_i,\ldots,j_{l})\in D_{i_1,\ldots,i_{l-1}}  } c_{j_1,\ldots,j_{l}} \le 1$,
  \item $0\le t_i < t_j <1$ and $t_{i} + c_i \le t_{j}$ for any $i,j \in \pi_1 D$, $i< j$,
  \item $0\le t_{i_1,\ldots,i_{l-1},i_l} < t_{i_1,\ldots,i_{l-1},j_l}<1$  and $t_{i_1,\ldots,i_{l-1},i_l} + c_{i_1,\ldots,i_{l-1},i_l} \le t_{i_1,\ldots,i_{l-1},j_l}<1$ for any $l=2,\ldots,d$ and $(i_1,\ldots,i_l), (i_1,\ldots,j_l) \in D_{i_1,\ldots,i_{l-1}}$, $i_l < j_l$, 
  \item let $i=\max \{i\in \pi_1 D\}$ then  $t_{i}+c_i \le 1$,
  \item for any $l=2,\ldots,d$ and $(i_1,\ldots,i_{l-1}) \in \pi_{l-1}(D)$, let $i_l=\max \{j_l \colon (i_1,\ldots,i_{l-1},j_l)\in D_{i_1,\ldots,i_{l-1}}\}$, 
then $t_{i_1,\ldots,i_{l-1},i_l}+c_{i_1,\ldots,i_{l-1},i_l} \le 1$.
\end{itemize}
This is essentially saying that our maps will map the unit cube into itself and the IFS will satisfy the open set condition. Now for any $\mathbf{i}\in D$ define the affine contraction $S_{\textbf{i}} \colon [0,1]^d \rightarrow [0,1]^d$ by
\[
S_{\textbf{i}}(x_1,\ldots, x_d) = (c_{i_1} x_1 + t_{i_1},\ldots, c_{i_1,\ldots,i_d}x_d+t_{i_1,\ldots,i_d}).
\]
The attractor of the IFS $\left\{S_{\textbf{i}} \right\}_{\textbf{i}\in D}$ is called a Lalley-Gatzouras sponge. This sponge satisfies a weak ordering condition only when, for some $l=2,\ldots, d$, $c_{i_1,\ldots,i_l}=c_{i_1,\ldots,i_{l-1}}$ for all $(i_1,\ldots, i_{l-1})\in \pi_{l-1}D$ and $(i_1,\ldots, i_{l})\in D_{i_1,\ldots, i_{l-1}}$. In this case we group together all the coordinates that cannot be ordered as before into $d^*$ many sets $\mathcal{J}_k= \prod_{i=a_1+\cdots + a_{k-1}+1}^{a_1+\cdots +a_k} \mathcal{I}_i$ where $c_{i_1,\ldots, i_{a_1+\cdots+a_{k-1}+1}}=\cdots = c_{i_1,\ldots, i_{a_1+\cdots+a_{k}}}$ for all $(i_1,\ldots, i_{a_1+\cdots+a_{k-1}+1})\in D_{i_1,\ldots, i_{a_1+\cdots+a_{k-1}}},\ldots,(i_1,\ldots, i_{a_1+\cdots+a_{k}})\in D_{i_1,\ldots, i_{a_1+\cdots+a_{k}-1}}$ and $(i_1,\ldots,i_{a_1+\cdots+a_{k-1}})\in \pi_{a_{1}+\cdots+a_{k-1}}D$; then our maps can be thought of as having the following form
\[
S_{\textbf{i}}(\textbf{x}_1,\ldots, \textbf{x}_{d^*}) = (c_{i_1} \textbf{x}_1 + \textbf{t}_{i_1},\ldots, c_{i_1,\ldots,i_{d^*}}\textbf{x}_{d^*}+\textbf{t}_{i_1,\ldots,i_{d^*}}).
\]
In an abuse of notation we redefine $D_{i_1,\ldots,i_{a_1+\cdots +a_{l-1}}}$ to be
\[
D_{i_1,\ldots,i_{a_1+\cdots +a_{l-1}}} = \left\{(j_1,\ldots,j_{a_1+\cdots+ a_{l}}) \in \pi_l(D)\colon (i_1,\ldots, i_{a_1+\cdots +a_{l-1}})=(j_1,\ldots,j_{a_1+\cdots +a_{l-1}}) \right\}
\] 
and 
\[
D_l=\pi_l D
\]
where $\pi_l$ is the projection of $D$ onto $\prod_{k=1}^l \mathcal{J}_k$.

Finally let $s$ be the dimension of $\pi_1 D$ (a self-similar set satisfying the open set condition) and for any $\textbf{i}=(i_1,\ldots,i_{d}) \in D$ and $l=1,\ldots,d^*-1$, let $s(i_1,\ldots,i_{a_1+\cdots+a_l})$ be the unique number satisfying 
\[
\sum_{(i_1,\ldots,i_{a_1+\cdots+a_l+a_{l+1}})\in D_{i_1,\ldots,i_{a_1+\cdots+a_l}}} c_{i_1,\ldots,i_{a_1+\cdots+a_l+a_{l+1}}}^{s(i_1,\ldots,i_{a_1+\cdots+a_l})}=1.
\]

\section{Results} \label{results}

In this paper we will provide an answer to question \cite[Question 4.2]{fraser-howroyd} by calculating the Assouad dimension of Bedford-McMullen and Lalley-Gatzouras sponges with a weak coordinate ordering condition.

\begin{thm} \label{BedfordMcMullen}
Let $K$ be a Bedford-McMullen sponge with a weak coordinate ordering condition. The Assouad dimension of $K$ is 
\begin{equation}\label{new}
\dim_{\text{\emph{A}}} K \ = \ \frac{\log N}{\log n_1^*} \ + \ \sum_{l=2}^{d^*} \frac{\displaystyle\log\max_{(i_1,\ldots, i_{a_1+\cdots+a_{l-1}})\in D_{l-1} } N(i_1, \ldots, i_{a_1+\cdots+a_{l-1}})}{\log n^*_l} 
\end{equation}
and the lower dimension of $K$ is
\[
\dim_{\text{\emph{L}}} K \ = \ \frac{\log N}{\log n_1^*} \ + \ \sum_{l=2}^{d^*} \frac{\displaystyle\log\min_{(i_1,\ldots, i_{a_1+\cdots+a_{l-1}})\in D_{l-1} } N(i_1, \ldots, i_{a_1+\cdots+a_{l-1}})}{\log n^*_l} 
\]
\end{thm}

The dimension of a Bedford-McMullen sponge whose coordinates are well ordered, calculated in \cite{fraser-howroyd}, is
\begin{equation}\label{old}
\dim_\text{\emph{A}} K = \frac{\log N}{\log n_l}+ \sum_{l=2}^{d} \frac{\displaystyle\log \max_{(i_1,\ldots,i_{l-1})\in D_{l-1} } N'(i_i,\ldots, i_{l-1})}{\log n_l},
\end{equation}
where $N'(i_i,\ldots, i_{l-1})$ is the number of ways of choosing the next digit of $(i_i,\ldots, i_{l-1})$ as a subset of the original $\mathcal{I}$. Formula (\ref{old}) is found by considering each coordinate separately as a collection of 1-dimensional Cantor sets (formed by the different columns). For each coordinate one finds the set with the greatest similarity dimension and then we add these maxima together to obtain the Assouad dimension of the whole.

When the coordinates of our sponge cannot be strictly ordered as in our setting, instead of always considering 1-dimensional sets we sometimes look at higher dimensional self-similar sets for which we can again find the greatest similarity dimension. The actual Assouad dimension of our set is then just the sum of the maximum similarity dimensions. 

Comparing formulas (\ref{new}) and (\ref{old}), there are clearly sponges, satisfying just a weak ordering condition, whose dimensions, predicted by the formulas, will coincide. But there are many others whose dimension, determined by (\ref{new}), is strictly smaller than that predicted by (\ref{new}). In fact, formula (\ref{old}) can give two different dimensions when a sponge satisfies only a weak ordering condition as one could swap two coordinates with equal contraction ratio, changing the $N'(i_1,\ldots i_l)$. 

Equality of the two formulas occurs when the largest multi-dimensional columns consist of only the maximal 1-dimensional columns that appear in (\ref{old}), that is 
\begin{align*}
\max_{(i_1,\ldots, i_{a_1+\cdots+a_{l-1}})\in D_{l-1} } N(i_1, \ldots, &i_{a_1+\cdots+a_{l-1}})= \\
&\prod_{k=1}^{a_l}\max_{(i_1,\ldots,i_{a_1+\cdots+a_{l-1}+k})\in D_{a_1+\cdots+a_{l-1}+k} } N'(i_1,\ldots, i_{a_1+\cdots+a_{l-1}+k}).
\end{align*}

As mentioned previously this dimension drop does not occur in the box and Hausdorff dimensions, or the Assouad dimension of carpets, and as such it would be interesting to try and find other, non-carpet examples where such an issue manifests. 

We motivate this difference by two examples. Figure \ref{fig:part-affine} is a sponge $F$ whose coordinates cannot be strictly ordered but whose Assouad dimension can be calculated using either formula (\ref{new}) or (\ref{old}):
\[
\dim_{\text{\emph{A}}} F = 1 + \frac{\log 3}{\log 3} = 2.
\]
However if we slightly modify that sponge by adding one extra map so $D= \left\{(0,0,0),(0,1,1),(0,2,1),(0,2,2),(1,0,1) \right\}$ with the same contraction ratios then we notice that the formulas for the Assouad dimension of $F$ provide two different values; the one in \cite{fraser-howroyd} gives
\[
2+ \frac{\log 2}{ \log 3} \approx 2.63 
\]
whereas the real value of the Assouad dimension is 
\[
\dim_{\text{\emph{A}}} F = 1 + \frac{\log 4}{\log 3} \approx 2.26.
\]
A third example of this phenomenon can be found in \cite[Section 4.3]{fraser-howroyd} where the Assouad dimension of a sponge with this weak coordinate ordering condition was calculated using the fact that the particular sponge was a product of two self-similar sets.

The Assouad and lower dimensions of Lalley-Gatzouras sponges were recently calculated in \cite{dfsu} but their proof relies on a condition analogous to the strict coordinate ordering needed in \cite{fraser-howroyd}. A simple extension of our result provides the following.

\begin{thm} \label{LalleyGatzouras}
Let $K$ be a Lalley-Gatzouras sponge with a weak coordinate ordering condition. The Assouad dimension of $K$ is 
\begin{equation}\label{lalley}
\dim_{\text{\emph{A}}} K \ = \ s \ + \ \sum_{l=2}^{d^*} \max_{(i_1,\ldots, i_{a_1+\cdots+a_{l-1}})\in D_{l-1} } s(i_1, \ldots, i_{a_1+\cdots+a_{l-1}})
\end{equation}
and the lower dimension of $K$ is
\[
\dim_{\text{\emph{L}}} K \ = \ s \ + \ \sum_{l=2}^{d^*} \min_{(i_1,\ldots, i_{a_1+\cdots+a_{l-1}})\in D_{l-1} } s(i_1, \ldots, i_{a_1+\cdots+a_{l-1}}).
\]
\end{thm}

Interestingly, Lalley and Gatzouras did not allow for any equality of contractions when computing the Hausdorff and box dimensions, however Das and Simmons \cite[Corolloary 3.4]{das-simmons} calculated the Hausdorff dimensions of a much more general set of sponges, which include Lalley-Gatzouras sponges that satisfy a weak ordering condition and no drop is perceived for these dimensions.

\section{Proof} \label{proof}

In this section we will prove Theorems \ref{BedfordMcMullen} and \ref{LalleyGatzouras}. In \ref{notation} we shall introduce any remaining notation that is needed for the proof of Theorem \ref{BedfordMcMullen}. Then we will prove the upper bound for the Assouad dimension of Bedford-McMullen sponges in \ref{upperbound} followed by the lower bound in \ref{lowerbound}. For completeness we will include a proof of both the upper and lower bounds for the lower dimension in \ref{lowerlower} and \ref{upperlower}. Much of this part will be similar to the proof in \cite{fraser-howroyd} so we will omit many of the technical details that have already been covered in that paper but we will nonetheless provide a summary of these techniques. Finally in Section \ref{LalleyGatzourasProof} we will sketch how the Bedford-McMullen proofs can be extended to the Lalley-Gatzouras case.

\subsection{Notation for Bedford-McMullen sponges with a weak coordinate ordering condition} \label{notation}

We define $\sigma:D^{\mathbb{N}} \rightarrow D^{\mathbb{N}}$ to be the shift map $\sigma(\textbf{i}_1,\textbf{i}_2,\ldots)=(\textbf{i}_2, \textbf{i}_3, \ldots)$, which acts in the inverse direction of $S_{\textbf{i}}$ but on the symbolic space instead of the unit cube.

Cylinder sets are unfortunately not optimal covers, as such we introduce the approximate cubes which are similar to cylinders and regular cubes at the same time, giving us the tools from the IFS with the symbolic space and the optimal covers from cubes. These will be used extensively throughout our proofs. For all $r\in (0,1]$ we choose the unique integers $k_1(r),\ldots,k_{d}(r)$, greater than or equal to 0, such that
\[
\left(\frac{1}{n_l}\right)^{k_l(r)+1}< r \leq \left(\frac{1}{n_l}\right)^{k_l(r)}
\]
for $l=1,\ldots,d$. In particular, $ \frac{-\log r}{\log n_l}-1 < k_l(r) \leq \frac{-\log r}{\log n_l}$. We observe that $k_{a_1+\cdots+a_{l-1}+1}=\cdots = k_{a_1+\cdots+a_l}$ and denote the common value by $k_l^*(r)=k_{a_1+\cdots+a_l}(r)$ for any $l=1,\ldots, d^*$; this is just the integer associated with all the coordinates in $\mathcal{J}_l$. Then the approximate cube $Q(\omega, r)$ of (approximate) side length $r<1$ determined by $\omega =\left( \textbf{i}_1, \textbf{i}_2 , \ldots \right) =\left( (i_{1,1}, \dots, i_{1,d}), (i_{2,1}, \dots, i_{2,d}) , \ldots \right)    \in D^{\mathbb{N}}$ is defined by
\[
Q(\omega, r)=\left\{ \omega'=\left( \textbf{j}_1, \textbf{j}_2 , \ldots \right)\in D^{\mathbb{N}} : \forall \, \,  l=1, \ldots, d \text{ and } \forall\, \, t= 1, \ldots, k_l(r) \text{ we have } j_{t,l}=i_{t,l} \right\}.
\]
The geometric analogue of approximate cubes, $\tau\left(Q(\omega, r)\right)$, is slightly harder to define and is contained in
\[ 
\prod_{l=1}^d \left[\frac{i_{1,l}}{n_l}+\cdots+\frac{i_{k_l(r),l}}{n_l^{k_l(r)}} \, , \, \frac{i_{1,l}}{n_l}+\cdots+\frac{i_{k_l(r),l}}{n_l^{k_l(r)}}+\frac{1}{n_l^{k_l(r)}} \right];
\]
 a rectangle in $\mathbb{R}^d$ aligned with the coordinate axes and of side lengths $n_l^{-k_l(r)}$, which are all comparable to $r$ since $ r \leq n_l^{-k_l(r)} < n_l r$. 

The lower bound for the Assouad dimension will use tangents and is interested in the `distance' between sets. This notion of distance is simply the Hausdorff metric $d_\mathcal{H}$ defined on the space of non-empty compact subsets of $\mathbb{R}^d$, which is defined by
\[
d_\mathcal{H}(A,B) \ = \  \inf \big\{ \varepsilon \geq 0   :   A \subseteq [B]_{\varepsilon} \text{ and } B \subseteq [A]_{\varepsilon} \big\}
\]
where $[A]_{\varepsilon}$ is the closed $\varepsilon$-neighbourhood of a set $A$.

\subsection{Upper bound for Assouad dimension} \label{upperbound}

Generally one uses a covering argument to obtain the upper bound, however our sponges being $d$-dimensional objects makes this idea quite difficult to implement. Instead we can simply use a measure theoretic technique thanks to the following proposition first proved in \cite[Proposition 3.1]{fraser-howroyd} and motivated by the measure theoretic definition developed in \cite{luksak, konyagin}:
\begin{prop}[\cite{fraser-howroyd}]\label{adupmeasure}
Suppose there exists a Borel probability measure $\nu$ on $D^{\mathbb{N}}$ and constants $C >0$ and $s \geq 0$ such that for any $0< r<R\le1$ and $\omega \in D^{\mathbb{N}}$ we have
\[
\frac{\nu \left( Q(\omega,R)\right)}{\nu \left( Q(\omega,r)\right)}\le C \left( \frac{R}{r} \right)^s.
\]
Then $\dim_{\text{\emph{A}}} K \le s$.
\end{prop}

The proof of this technical lemma is simple and follows from the original definition of Assouad dimension; for its proof see the original paper. Essentially it says that if one has a measure satisfying some condition on approximate cubes then we have an upper bound on the Assouad dimension. Thankfully we can modify the measure used to this end in \cite{fraser-howroyd} such that our desired dimension is an upper bound.

The modified measure is defined by
\[
p_\textbf{i}=p_{(i_1,\ldots,i_d)}=\frac{1}{N\left(\prod_{l=2}^{d^*} N(i_1,\ldots,i_{a_1+\cdots+a_{l-1}})\right)}.
\]
Again, this is simply the `coordinate uniform measure' introduced for the regular sponges except for the coordinates with equal contractions where it is a sort of uniform distribution, as desired. Similarly one can check from the definitions that $\sum_{ \textbf{i} \in D} p_\textbf{i} = 1$. We call this measure the `partial coordinate uniform measure', since it is a modified version of the original.

There is a precise formula for conditional probabilities on the symbol space as noted in Olsen \cite[Section 3.1]{sponges}, we omit the exact formulation as it simply reduces to, whenever $(i_1, \ldots, i_{a_1+\cdots + a_{l}})\in D_l$, $p(i_{a_1+\cdots + a_{l-1}+1},\ldots,i_{a_1+\cdots + a_{l}}  \vert i_{1},\ldots,i_{a_1+\cdots+a_{l-1}})=1/N(i_1,\ldots,i_{a_1+\cdots+a_{l-1}})$ for $l=2,\ldots,d^*$ and $p( i_1, \ldots, i_{a_1} \vert \emptyset)=1/N$. This conditional probability is just the probability of picking the next set of digits in the set $\mathcal{J}_l$ given the previous coordinates.

Technically our measure is defined on the pre-fractals of our sponge, not on approximate cubes. Fortunately Olsen \cite[(6.2)]{sponges} noted the following straightforward formula
\begin{equation} \label{approxcubemeasure}
\tilde \mu(Q(\omega,r))=\prod^d_{l=1} \prod_{j=0}^{k_l(r)-1}p_l(\sigma^j\omega)
\end{equation}
where $p_l(\omega)=p(i_{1,l}\vert i_{1,1},\ldots,i_{1,l-1})$. The important fact to note here is that $p(i_{l}\vert i_{1},\ldots,i_{l-1}) \times p(i_{l+1}\vert i_{1},\ldots,i_{l})=p(i_l,i_{l+1}\vert i_{1},\ldots,i_{l-1})$ and as such we can combine the conditional probabilities in (\ref{approxcubemeasure}) when two successive coordinates have the same contraction ratio so $k_l(r)=k_{l+1}(r)$.  We are now ready to obtain the upper bound.

\begin{proof} By (\ref{approxcubemeasure}), when $s$ is the desired Assouad dimension, we have
\begin{align*}
\frac{\tilde \mu(Q(\omega,R))}{\tilde \mu(Q(\omega,r))}
&=\frac{\prod^d_{l=1} \prod_{j=0}^{k_l(R)-1}p_l(\sigma^j\omega)}{\prod^d_{l=1} \prod_{j=0}^{k_l(r)-1}p_l(\sigma^j\omega)}  \\
&=\prod^d_{l=1} \prod_{j=k_l(R)}^{k_l(r)-1}\frac{1}{p_l(\sigma^j\omega)} \\
&=\prod^{d^*}_{l=1} \prod_{j=k_l^*(R)}^{k_l^*(r)-1}\frac{1}{p(i_{j+1,a_1+\cdots+a_{l-1}+1},\ldots,i_{j+1,a_1+\cdots+a_{l}}\vert i_{j+1,1},\ldots, i_{j+1,a_1+\cdots+a_{l-1}})} \\
&\leq \left( \prod_{j=k_1^*(R)}^{k_1^*(r)-1} N \right)\left(\prod^{d^*}_{l=2} \prod_{j=k_l^*(R)}^{k_l^*(r)-1} \max_{(i_1,\ldots, i_{a_1+\cdots+a_{l-1}})\in D_{l-1} } N(i_1, \ldots, i_{a_1+\cdots+a_{l-1}}) \right) \\
&=   N^{k_1^*(r)-k_1^*(R)} \left(\prod^{d^*}_{l=2} \max_{(i_1,\ldots,i_{l-1})\in D_{l-1}} N(i_1,\ldots,i_{l-1})^{k_l^*(r)-k_l^*(R)}\right) \\
& \le N^{\log R/\log n_1^*-\log r/\log n_1^*+1} \\
&\qquad \left(\prod^{d^*}_{l=2} \max_{(i_1,\ldots, i_{a_1+\cdots+a_{l-1}})\in D_{l-1} } N(i_1, \ldots, i_{a_1+\cdots+a_{l-1}})^{\log R/\log n_l^*-\log r/\log n_l^*+1}\right) \\
& \le \left(n^*_1\right)^{a_1} \times \cdots \times \left(n^*_{d^*}\right)^{a_d}\left(\frac{R}{r}\right)^{\displaystyle\frac{\log N}{\log n_1^*}}
\\
& \qquad \left(\prod^{d^*}_{l=2} \left(\frac{R}{r}\right)^{\frac{\displaystyle\log \max_{(i_1,\ldots, i_{a_1+\cdots+a_{l-1}})\in D_{l-1} } N(i_1, \ldots, i_{a_1+\cdots+a_{l-1}})}{\displaystyle\log n_l^*}}\right) \\
&\le n_d^{d}\left( \frac{R}{r}\right)^{s}.
\end{align*}

This estimate combined with Proposition \ref{adupmeasure} gives us the desired upper bound.

\end{proof}

\subsection{Lower bound for Assouad dimension} \label{lowerbound}

For the lower bound we will use `weak tangents', a technique possible thanks to a proposition of Mackay and Tyson \cite[Proposition 6.1.5]{mackaytyson} which was then modified by Fraser \cite[Proposition 7.7]{fraser}. This is the section where the original proof in \cite{fraser-howroyd} does not apply, there is however a way around this using some of the ideas detailed so far.

\begin{prop}[Very weak tangents]\label{tangents}
Let $X\subset \mathbb{R}^d$ be compact and let $F$ be a compact subset of $X$. Let $(T_k)$ be a sequence of bi-Lipschitz maps defined on $\mathbb{R}^d$ with Lipschitz constants $a_k, b_k \ge 1$ such that 
\[
a_k \lvert x-y \rvert \leq \lvert T_k(x) - T_k(y) \rvert \leq b_k \lvert x-y \rvert \,\,\,\,\,\,\,\, (x,y\in\mathbb{R}^d)
\]
and
\[
\sup_k b_k / a_k = C_0 <\infty
\]
and suppose that $T_k(F) \cap X \rightarrow \hat{F}$ in the Hausdorff metric. Then the set $\hat F$ is called a \emph{very weak tangent} to $F$ and, moreover,  $\dim_{\text{\emph{A}}} F \geq \dim_{\text{\emph{A}}} \hat{F}$.
\end{prop}

To simplify notation, we choose $\textbf{i}(l) = (i(l)_1, \ldots, i(l)_d)\in D$ for $l=2, \ldots d^*$ to be an element of $D$ which attains the maximum value for $N(i_1, \ldots, i_{a_1+\cdots+a_{l-1}})$, i.e.
\[
N \left(i(l)_1, \ldots, i(l)_{a_1+\cdots+a_{l-1}} \right)  \ =  \  \max_{(i_1,\ldots,i_{a_1+\cdots+a_{l-1}})\in D_{l-1}} N(i_1, \ldots, i_{a_1+\cdots+a_{l-1}}).
\]
There might be several possibilities for each $\textbf{i}(l)$, but we can pick it arbitrarily. 

The rest of this part will be a brief explanation of the technique as it is very similar to that in \cite{fraser-howroyd}. The main idea is to construct a tangent $\hat{F}$ that is some product of self-similar sets whose box dimension can be calculated using common formulas as in \cite{falconer} and the box dimension of $\hat{F}$ is simply the desired Assouad dimension. To show that this tangent is indeed a tangent one shows that $T_k(F) \cap X$ and $\hat{F}$ are both subsets of a set that is the product of some Menger-like sponge and $d^*-1$ pre-fractals of Menger-like sponges which depend on their respective coordinate collection (and the level depends on $k$) and then as $k$ tends to infinity the Hausdorff distance between all these sets becomes 0. There is also a complication that might arise due to the definition of approximate cubes but this is easily circumvented in the Assouad dimension case.

The desired tangent is 
\[
\hat{K}=\pi_1 K \times \prod_{l=2}^{d^*} K_l
\]
where $K_l$ is the Menger-like sponge obtained by the IFS acting on $[0,1]^{a_l}$ which forms a sponge identical to the maximal $a_l$-dimensional plane in the $a_1+\cdots + a_{l-1}+1$ to $a_1+\cdots + a_{l}$ coordinates of our original sponge. These are simply the Cantor sets of the original proof when the contraction ratios are different but self-similar sponges when they are equal. The set $\pi_1 K$ is the geometric projection of the sponge on to the first $a_1$ coordinates (note that we use the projection function $\pi$ in both the geometric and symbolic spaces). 

For $l=2, \dots, d^*$ and $m \in \mathbb{N}$, we let $K_l^m$ be the $m^{\text{th}}$ pre-fractal of $K_l$ where our initial set is $[0,1]^{a_l}$. In particular, the set $K_l^m$ is a union of $N(i_1, \ldots, i_{a_1+\cdots+a_{l-1}})^m$ cubes of side lengths $(n_l^*)^{-m}$.

Given a geometric approximate cube $\tau(Q) =\tau(Q(\omega,r))$, we define a bi-Lipschitz map $T^Q: \tau(Q)  \to [0,1]^d$ (or $T^Q:  \mathbb{R}^d  \to \mathbb{R}^d$) by
\[
T^Q(\textbf{x})= \begin{pmatrix}n_1^{k_1(r)}\left( x_1-\left( \frac{i_{1,1}}{n_1}+\ldots+\frac{i_{k_1(r),1}}{n_1^{k_1(r)}} \right) \right)\\\vdots\\n_d^{k_d(r)}\left( x_d-\left( \frac{i_{1,d}}{n_d}+\ldots+\frac{i_{k_d(r),d}}{n_d^{k_d(r)}} \right) \right)\end{pmatrix}.
\]
Thus $T^Q$ translates $\tau(Q)$ such that the point closest to the origin from the rectangle containing $\tau(Q)$ becomes the origin and then scales it up by a factor of $n_l^{k_l(r)}$ in each coordinate $l$. Thus these maps take the natural rectangle containing $\tau(Q)$ precisely to the unit cube $[0,1]^d$. These maps clearly satisfy the conditions imposed by Proposition \ref{tangents}, i.e., they are restrictions of bi-Lipschitz maps on $\mathbb{R}^d$ with constants $b_Q=\sup_{l=1, \ldots, d}n_l^{k_l(r)}$ and $a_Q=\inf_{l=1, \ldots, d}n_l^{k_l(r)}$ satisfying
\[
\frac{b_Q}{a_Q}\leq \sup_{l=1,\ldots,d}\frac{rn_l}{r}\leq n_d < \infty
\]
for any $Q$. This follows from the definition of $k_l(r)$ and does not rely on the inequality of the $n_l$.

We define, for small $R$, $\omega(R)=\left(\textbf{i}_1, \textbf{i}_2, \ldots \right) \in D^{\mathbb{N}}$ where $\textbf{i}_t = (i_{t,1}, \dots, i_{t,d}) =\textbf{i}(l)$ for $t=k_l^*(R)+1, \ldots,k_{l-1}^*(R)$ for all $l=2, \ldots, d^*$. So $\omega(R)$ has the form
\[
\omega(R)=\left( \textbf{i}_1, \ldots, \underbrace{\textbf{i}(d), \ldots, \textbf{i}(d)}_{k_{d^*-1}^*(R)-k_{d^*}^*(R) \text{ times}}, \ \ \underbrace{\textbf{i}(d-1), \ldots, \textbf{i}(d-1)}_{k_{d^*-2}^*(R)-k_{d^*-1}^*(R) \text{ times}}, \ldots,  \underbrace{\textbf{i}(2), \ldots, \textbf{i}(2)}_{k_1^*(R)-k_{2}^*(R) \text{ times}}, \ldots  \right).
\]

We can then prove the following technical lemma, see \cite[Lemma 3.3]{fraser-howroyd} for details.

\begin{lma}[\cite{fraser-howroyd}] \label{productform}
For $ R \in (0,1]$ small enough and $Q=Q(\omega(R),R)$, we have
\[
T^Q(\tau(Q)) \ \subseteq \  \pi_1 K \times \prod_{l=2}^{d^*} K_l^{k_{l-1}^*(R)-k_{l}^*(R)}.
\]
\end{lma}

By using this lemma and the definition of pre-fractals we can show that 
\[
d_{\mathcal{H}}\left( \hat{K},  \ T^Q(\tau(Q(\omega(R),R))) \right)\rightarrow 0
\]
as $R \to 0$.

One final technical lemma then solves the following: 
\[
T^Q(K) \cap [0,1]^d \ \supseteq \ T^Q(\tau(Q(\omega(R),R))) \ \to \ \hat{K},
\]
where the containment may be strict, which implies that $\hat{K}$ is not necessarily a tangent of $K$. 

\begin{lma}
Let $E_k, \,  F_k \subseteq [0,1]^d$ be  sequences of non-empty compact sets which converge in the Hausdorff metric to compact sets $E$ and $F$ respectively.   If $F_k \subseteq E_k$ for all $k$, then $F \subseteq E$.
\end{lma}

The proof of this classical lemma is left to the reader. We can now finally calculate the lower bound.

\begin{proof} Standard results on the box dimensions of product sets \cite[Chapter 7]{falconer} and of self-similar sets \cite[Chapter 9]{falconer} imply that
\begin{eqnarray*}
\dim_{\text{B}} \hat{K} &=& \dim_{\text{B}} \pi_1K+\sum_{l=2}^{d^*} \dim_{\text{B}} K_l \\ \\&=& \frac{\log N}{\log n_1}+\sum_{l=2}^{d^*} \frac{\displaystyle\log\max_{(i_1,\ldots, i_{a_1+\cdots+a_{l-1}})\in D_{l-1} } N(i_1, \ldots, i_{a_1+\cdots+a_{l-1}})}{\log n^*_l}.
\end{eqnarray*}

By compactness we can take a convergent subsequence of $T^Q(K) \cap [0,1]^d \rightarrow K'$. Then by Proposition \ref{tangents} there is a set $\hat{K} \subseteq K'$ and by monotonicity of Assouad dimension the following holds
\begin{eqnarray*}
\dim_{\text{A}} K  \  \ge \ \dim_{\text{A}} K'  \  \ge \ \dim_{\text{A}} \hat{K}  \ge &\dim_{\text{B}} \hat{K}, 
\end{eqnarray*}
as required.
\end{proof}

\subsection{Lower bound for lower dimension} \label{lowerlower}

The proof of the lower bound for the lower dimension is essentially the same as the proof in \ref{upperbound}. We start by introducing a lemma proved in \cite[Proposition 3.5]{fraser-howroyd}.

\begin{prop}[\cite{fraser-howroyd}]\label{lowlowmeasure}
Suppose there exists a Borel probability measure $\nu$ on $D^{\mathbb{N}}$ and constants $C >0$ and $s \geq 0$ such that for any $0< r<R\le1$ and $\omega \in D^{\mathbb{N}}$ we have
\[
\frac{\nu \left( Q(\omega,R)\right)}{\nu \left( Q(\omega,r)\right)}\ge C \left( \frac{R}{r} \right)^s.
\]
Then $\dim_{\text{\emph{L}}} K \ge s$.
\end{prop}

Then using a mixture of the original lower dimension's lower bound proof and the new Assouad dimension upper bound proof we get

\begin{proof}
\begin{align*}
\frac{ \tilde \mu(Q(\omega,R))}{\tilde \mu(Q(\omega,r))}
&\ge  n_d^{-d}\left( \frac{R}{r}\right)^{\displaystyle\frac{\log N}{\log n_1^*} \ + \ \sum_{l=2}^{d^*} \frac{\displaystyle\log\min_{(i_1,\ldots, i_{a_1+\cdots+a_{l-1}})\in D_{l-1} } N(i_1, \ldots, i_{a_1+\cdots+a_{l-1}})}{\log n^*_l} }.
\end{align*}

This estimate combined with Proposition \ref{lowlowmeasure} gives us the required lower bound.
\end{proof}

\subsection{Upper bound for lower dimension} \label{upperlower}

The method for this final proof is again a mixture of the original one with the new idea to avoid the equality problem and should be clear from these other proofs, as such we leave it to the reader.

\subsection{Extension to Lalley-Gatzouras sponges}\label{LalleyGatzourasProof}

In this section we will introduce the more generalised notation and concepts needed to calculate the Assouad dimension of Lalley-Gatzouras sponges. The proof for the lower dimension then follows easily and as such is skipped.

The definition of $k_l(r)$ cannot be simply copied to Lalley-Gatzouras sponges as the contractions for any given coordinate are not necessarily uniform, as such, given a word $\omega=(\textbf{i}_1,\textbf{i}_2,\ldots) \in D^{\mathbb{N}}$, $r\in (0,\min_{(i_1,\ldots,i_d)\in D} c_{i_1,\ldots,i_d}]$ and $l=1,\ldots,d $ we choose $k_l(r,\omega)$ to be the unique integer, greater than or equal to 1, such that
\[
\prod_{m=1}^{k_l(r,\omega)+1} c_{i_{m,1},\ldots,i_{m,l}} < r \le \prod_{m=1}^{k_l(r,\omega)} c_{i_{m,1},\ldots,i_{m,l}}.
\]
Note that when two coordinates, say $l$ and $l+1$, have at least one $\textbf{j}\in D$ such that $c_{j_1,\ldots,j_l} > c_{j_1,\ldots,j_{l+1}}$ then there exists a word $\omega=(\textbf{j},\ldots)$ such that $k_l(r,\omega)>k_{l+1}(r,\omega)$. When no such contraction exists then our sponge satisfies a weak ordering condition and we identify the two coordinates as one as previously explained. As before we identify $k^*_l(r,\omega)=k_{a_1+\cdots+a_l}(r,\omega)$.

The definition of approximate cubes is the same as for Bedford-McMullen sponges; note that the geometric analogue of approximate cubes for Lalley-Gatzouras sponges are contained in a slightly more notationally complex product of intervals, but the general concept still holds. 

To calculate the upper bound for the Assouad dimension we again use the measure theoretic definition using the following measure
\[
p_\mathbf{i}=p_{i_1,\ldots,i_d}=c_{i_1,\ldots,i_{a_1}}^s c_{i_1,\ldots,i_{a_1+a_2}}^{s(i_1,\ldots,i_{a_1})}\cdots c_{i_1,\ldots ,i_{d}}^{s(i_1,\ldots,i_{a_1+\cdots+a_{d^*-1}})}.
\]
We call this the `coordinate measure of full dimension', as, like the coordinate uniform measure, it distributes the mass in such a way that the measure `sees' the dimension of each coordinate individually, up to the weak ordering, where the dimension of the grouped coordinates is observed. 

Naturally the conditional probabilities are what we expect: 
\[
p(i_{a_1+\cdots+a_{l-1}+1},\ldots,i_{a_1+\cdots+a_l}\colon i_1,\ldots,i_{a_1+\cdots+a_{l-1}})=c_{i_1,\ldots,i_{a_1+\cdots+a_l}}^{s(i_1,\ldots,i_{a_1+\cdots+a_{l-1}})}.
\]
Thus the following calculation provides our upper bound.
\begin{align*}
\frac{\tilde{\mu}(Q(\omega,R))}{\tilde{\mu}(Q(\omega,r))}& = \frac{\prod_{l=1}^d\prod_{j=0}^{k_l(R,\omega)-1}p_l(\sigma^j\omega)}{\prod_{l=1}^d\prod_{j=0}^{k_l(r,\omega)-1}p_l(\sigma^j\omega)} \\ 
&=\left(\prod_{j=k_1^*(R,\omega)}^{k_1^*(r,\omega)-1}c_{i_{j,1},\ldots,i_{j,a_1}}^{-s} \right)\times\left( \prod_{l=2}^{d^*}\prod_{j=k_l^*(R,\omega)}^{k_l^*(r,\omega)-1}c_{i_{j,1},\ldots,i_{j,a_1+\cdots+a_l}}^{-s(i_{j,1},\ldots,i_{j,a_1+\cdots+a_{l-1}})}\right)\\ 
& \le \left(\prod_{j=k_1^*(R,\omega)}^{k_1^*(r,\omega)-1}c_{i_{j,1},\ldots,i_{j,a_1}}\right)^{-s} \times \prod_{l=2}^{d^*}\left( \prod_{j=k_l^*(R,\omega)}^{k_l^*(r,\omega)-1}c_{i_{j,1},\ldots,i_{j,a_1+\cdots+a_l}}\right)^{-\max s(i_{j,1},\ldots,i_{j,a_1+\cdots+a_{l-1}})}\\
& \le \left( \min_{\mathbf{i}\in D} c_{i_1,\ldots,i_d}^{-d} \right) \left( \frac{R}{r}\right)^s\times \prod_{l=2}^{d^*}\left( \frac{R}{r}\right)^{\max s(i_{j,1},\ldots,i_{j,a_1+\cdots+a_{l-1}})}
\end{align*}
where the exponent is our desired dimension. Technically to use the measure theoretic definition we need an analgous version of Proposition \ref{adupmeasure} which was originally proved just for Bedford-McMullen sponges, however the proof for Lalley-Gatzouras sponges is the same and is omitted.
 
For the lower bound we wish to use Proposition \ref{tangents}. The technique for this part is the same as Section \ref{lowerbound} in the sense that we want to construct a weak tangent that is a product of self-similar sets. We define $T^Q$, where $Q$ is an approximate cube, to be the function which maps $Q$ to the unit cube as before. 

For each $l=2,\ldots,d^*$, define $\mathbf{i}(l)=(i(l)_1,\ldots,i(l)_d)\in D$ to be an element such that 
\[
s(i(l)_1,\ldots,i(l)_{a_1+\cdots+a_{l-1}})=\max_{(i_1,\ldots,i_{a_1+\cdots+a_{l-1}})\in D_{l-1}}s(i_1,\ldots,i_{a_1+\cdots+a_{l-1}}).
\]
We also let $\mathbf{i}(\text{twist},l)=(i(\text{twist},l)_1,\ldots,i(\text{twist},l)_d)\in D$ be an element of $D$ such that 
\[
c_{i(\text{twist},l)_1,\ldots,i(\text{twist},l)_{a_1+\ldots+a_{l-1}}}> c_{i(\text{twist},l)_1,\ldots,i(\text{twist},l)_{a_1+\ldots+a_l}}.
\]
This condition will be used to make $k_l^*(R,\omega)<k_{l-1}^*(R,\omega)$ and is possible due to our identification of coordinates.

Heuristically the approximate cubes $Q(\omega(R),R)$ that we zoom into for our weak tangent will be the same as before, that is from $k_l^*(R,\omega(R))$ to $k_{l-1}^*(R,\omega(R))$ we want to pick the largest column for the $l^{\text{th}}$ coordinate. However we want the $k_{l-1}^*(R,\omega(R))-k_l^*(R,\omega(R))$ to be positive and tend to infinity as $R$ goes to zero; there exists Lalley-Gatzouras sponges where this does not always happen without an additional step. To overcome this we choose the first $k_{d}^*(R,\omega(R))$ elements of $\omega(R)$ to be our twist symbols, so as $R$ goes to 0 the differences will tend to infinity.

Formally, for small enough $R$, we define $\omega(R)=(\mathbf{i}_1,\mathbf{i}_2,\ldots)\in D^\mathbb{N}$ by setting 
\[
\mathbf{i}_{m d^* + 1},\ldots,\mathbf{i}_{m d^* + d^*}= \mathbf{i}(\text{twist},1), \ldots,\mathbf{i}(\text{twist},d^*)
\]
for $m=0,1,\ldots, \lfloor\frac{k_d^*(R,\omega(R))}{d^*}\rfloor-1$ and $\mathbf{i}_t=(i_{t,1},\ldots,i_{t,d}) = \mathbf{i}(l)$ for all $t=k^*_l(R,\omega(R))+1,\ldots, k^*_{l-1}(R,\omega(R))$ and for all $l=2,\ldots, d^*$. Due to the definition of $k_l^*(R,\omega)$, fixing elements after the $k_l^*$-th symbol does not change $k_l^*(R,\omega)$, so such a word is well-defined.

To finish, in much the same way as in Section \ref{lowerbound}, one can check that $T^{Q(R,\omega(R))}(K) \cap [0,1]^d$ will converge to a set with a subset of the desired dimension.

\vspace{10pt}
\begin{centering}

\textbf{Acknowledgements}

\end{centering}

The author was supported by the EPSRC Doctoral Training Grant EP/N509759/1 whilst writing this paper. He thanks his supervisor Jonathan Fraser for the advice received during the project and an anonymous referee for many stimulating comments.

\end{document}